\documentclass[11pt,x11names]{article}
\usepackage{lmodern}
\usepackage{amssymb,bbm}    \usepackage{amsmath}    \usepackage{amsthm} \usepackage{mathrsfs}
\usepackage{tikz,color}  \usepackage{mathtools}
\usepackage{enumitem}
\usepackage{bookmark}
\usepackage{array} \usepackage[font=footnotesize]{caption}
\usepackage{colortbl}
\usepackage{subcaption}
\usepackage{placeins}
\newtheorem{theorem}{Theorem}[section]

\newtheorem{lemma}[theorem]{Lemma}

\theoremstyle{remark}
\newtheorem{remark}[theorem]{Remark}

\usepackage{hyperref}
\hypersetup{colorlinks=true,urlcolor=black,linkcolor=black,citecolor=black,naturalnames=true,hypertexnames=false}

\usepackage[backend=bibtex,giveninits, style=authoryear, uniquelist=false, natbib=true, doi=false, url=false, isbn = false,maxbibnames=99]{biblatex}
\DeclareFieldFormat[article]{title}{#1} \DeclareFieldFormat[book]{title}{\itshape #1}
\DeclareFieldFormat[thesis]{title}{#1}
\DeclareFieldFormat[incollection]{title}{\itshape #1}
\DeclareFieldFormat[inproceedings]{title}{\itshape #1}
\DeclareFieldFormat[article]{author}{\color{blue} #1}
\DeclareFieldFormat[article]{pages}{#1.}
\DeclareFieldFormat[incollection]{pages}{#1.} 

\renewbibmacro{in:}{\ifentrytype{article}{}{\printtext{\bibstring{in}\intitlepunct}}}

\DeclareNameAlias{sortname}{last-first}
\DeclareNameAlias{default}{last-first}

\renewbibmacro*{volume+number+eid}{\printfield{volume}\setunit*{\addnbthinspace}\printfield{number}\setunit{\addcomma\space}\printfield{eid}}
\DeclareFieldFormat[article]{number}{\mkbibparens{#1}}

\renewbibmacro*{publisher+location+date}{\printlist{publisher}\iflistundef{location}
	{\setunit*{\addcomma\space}}
	{\setunit*{\addcomma\space}}\printlist{location}\setunit*{\addcomma\space}\usebibmacro{date}\newunit}

 \DeclareFieldFormat[article]{number}{}

\def\D{{\cal D}}

\def\R{{\mathbb R}}
\def\Rd{{\mathbb R}^d}

\def\PROB {{\mathbb P}}
\def\EXP {{\mathbb E}}
\def\IND{{\mathbb I}}
\def\Var{{\mathrm{Var}}}

\def \p {^{\prime}}

\newcommand{\ftilde}{{\widetilde f}}
\newcommand{\Fc}{{\mathcal F}}
\newcommand{\Qc}{{\mathcal Q}}
\newcommand{\KL}{{\mathrm{KL}}} \newcommand{\TV}{{\mathrm{TV}}}  \newcommand{\ib}{{\boldsymbol i}}
\newcommand{\jb}{{\boldsymbol j}}
\newcommand{\yb}{{\boldsymbol y}}
\newcommand{\xb}{{\boldsymbol x}}
\newcommand{\zb}{{\boldsymbol z}}
\newcommand{\Zb}{{\boldsymbol Z}}
\newcommand{\Flipd}{\Fc_{\mathrm{Lip}}^d}
\newcommand{\dd}{\mathrm{d}}
\renewcommand{\epsilon}{\varepsilon}
\newcommand{\mutilde}{{\widetilde \mu}}
\newcommand{\defeq}{\vcentcolon=}

\newcommand{\Zc}{{\mathcal Z}}
\newcommand{\Zs}{{\mathscr Z}}

\allowdisplaybreaks

\begin{document}

\title{Multivariate density estimation from privatised data: universal consistency and minimax rates}
\author{L\'aszl\'o Gy\"orfi\thanks{Department of Computer Science and Information Theory, Budapest University of Technology and Economics, Magyar Tud\'{o}sok krt. 2., Budapest, H-1117, Hungary.
\texttt{gyorfi@cs.bme.hu}}
\and
Martin Kroll\thanks{Ruhr-Universität Bochum, Fakultät für Mathematik, Universitätsstra\ss e 150. D-44801 Bochum.} $^{,}$\thanks{Universität Bayreuth, Fakultät für Mathematik, Physik und Informatik. D-95440 Bayreuth. \texttt{martin.kroll@uni-bayreuth.de}}}
\thispagestyle{empty}
\maketitle

\begin{abstract}
	We revisit the classical problem of nonparametric density estimation but impose local differential privacy constraints.
	Under such constraints, the original multivariate data $X_1,\ldots,X_n \in \mathbb{R}^d $ cannot be directly observed, and all estimators are functions of the randomised output of a suitable privacy mechanism.
	The statistician is free to choose the form of the privacy mechanism, and in this work we propose to add Laplace distributed noise to a discretisation of the location of an observed vector.
	Based on these randomised data, we propose a novel estimator of the density function, which can be viewed as a privatised version of the well-studied histogram density estimator.
	Our theoretical results include universal pointwise consistency and strong universal $L_1$-consistency.
	In addition, a convergence rate for Lipschitz continuous functions is derived, which is complemented by a matching minimax lower bound.
	We illustrate the trade-off between data utility and privacy by means of a small simulation study.
\end{abstract}


\smallskip


\normalsize

\section{Introduction}\label{s:intro}

Let $X$ be a  random  vector taking values in $\R^d$.
We denote by $ \mu $ the distribution of the vector $X$, that is, for all Borel sets $A\subset \R^d$, we have $\mu(A)=\PROB(X\in A)$.
Assume that $ \mu $ has a density $f$ with respect to the Lebesgue measure.
Classical nonparametric theory considers the problem of estimating $f$ based on data consisting of independent, identically distributed (i.i.d.) copies of the vector $X$,
\begin{align}\label{eq:D_n}
	\D_n=\{X_1,\ldots,X_n\}.
\end{align}
Several loss functions are typically considered in this context in order to evaluate the performance of an estimator $\check f_n$.
In this work, we consider the $L_1$-error given through
\begin{equation*}
	\int_{\R^d} |f(\xb) - \check f_n(\xb) | \dd \xb.
\end{equation*}
This choice is natural since densities belong by their very definition to the space $L_1(\R^d)$ of functions with finite $L_1$-norm.
We refer the interested reader to \citep{DeGy85} for a comprehensive introduction to density estimation from the $L_1$-point of view.

Here, we study density estimation under the $L_1$-error as well but we differ from the classical framework in the sense that we do no longer assume that the raw data \eqref{eq:D_n} are accessible.
Instead, only anonymised surrogates of the raw data which are supposed to satisfy the specific privacy condition of \emph{local differential privacy} (LDP) can be used for inference.

The main motivation of the present paper was the question whether \emph{strong universal consistency}, that is, almost sure convergence of the $L_1$-error to zero for any density function can be established for an estimator defined in terms of such anonymised data.
We answer this question in the affirmative using a histogram type estimator adapted to the privacy setup which will be introduced in Section~\ref{s:hist}.

In order to pigeonhole the methods and results obtained in the rest of the paper, let us shortly recap the definition and theoretical properties of histograms based on cubic partitions when the raw data \eqref{eq:D_n} are available.
Given a partition $\{A_{h_n,1},A_{h_n,2},\ldots\}$
of $\Rd$ consisting of cubes $A_{h_n,j}$ of volume $h_n^d$, the corresponding classical histogram estimator is defined by
\begin{equation}
	\label{hist}
	f_{n}(\xb)= \frac{\mu_n(A_{h_n,j})}
	{h_n^d}\qquad \mbox{if } \xb\in A_{h_n,j},
\end{equation}
where for any measurable set $A$ we denote by
\begin{align*}
	\mu_n(A) = \frac{1}{n} \sum_{i=1}^n \IND_{\{X_i\in A\}}
\end{align*}
the empirical counterpart of $\mu(A)$ and $\IND$ is the indicator function.
The following well-known theorem states that the estimator $f_n$ is strongly universally consistent.
\begin{theorem}[\citep{Abo76}; Theorem~3.2 in \citep{DeGy85}]
	\label{hcons}
	If
	\begin{align*}
		\lim_{n\to\infty}{ h}_n=0\, \mbox{ and }\, \lim_{n\to\infty}n{ h}_n^{d} =\infty,
	\end{align*}
	then
	\begin{equation*}
\lim_{n\to \infty} \int_{\R^d} |f(\xb) - f_n(\xb) | \dd \xb=0 \quad \text{almost surely.}
	\end{equation*}
\end{theorem}

Incidentally, we will also derive a minimax rate of convergence of our modified histogram estimator for Lipschitz continuous densities.
In the classical setting with accessible raw data $\mathcal D_n$, the following result (stated here in a slightly more general form than in the given references) is well-known.
\begin{theorem}[\citep{BeGy98}; Theorems 5.5 and 5.6 in \citep{DeGy85}]
	\label{hrate}
	If $X$ has a bounded support $S$ of Lebesgue measure $\lambda(S)$, and $f$ is Lipschitz continuous with Lipschitz constant $L$, then
	\begin{align*}
		\EXP \left[  \int_{\R^d} |f(\xb) - f_n(\xb) | \dd \xb \right] 
		\le
		L\sqrt{d}h_n+\frac{\lambda(S)}{\sqrt{nh_n^d}}.
	\end{align*}
	In particular, the choice $h_n=c\cdot n^{-1/(d+2)}$ yields
	\begin{align*}
		\EXP\left[  \int_{\R^d} |f(\xb) - f_n(\xb) | \dd \xb \right] = O\left(n^{-1/(d+2)}\right).
	\end{align*}
\end{theorem}
In addition, the rate $n^{-1/(d+2)}$ can be identified as the optimal rate of convergence for Lipschitz continuous functions.

The problem of density estimation under privacy restrictions has already been tackled by statisticians before since the notion of differential privacy was introduced in \citep{dwork2006differential} so that a short review of relevant literature seems necessary in order to fit our work into the already existing ones. The first paper dealing with density estimation under differential privacy is \citep{wasserman2010statistical}. The authors of that paper, however, worked in the setup of \emph{global} differential privacy and under the less restrictive condition of \emph{approximate} differential privacy. The work on nonparametric density estimation under \emph{local} differential privacy has been initiated afterwards in the important paper \citep{DuJoWa18} where minimax rates of convergence of convergence over Sobolev ellipsoids were established. These results were extended to more general ellipsoids in Besov spaces in \citep{BuDuKrSa20} where wavelet estimators under LDP were studied. The estimation of linear functionals under local differential privacy was investigated in great detail in \citep{rohde2020geometrizing}.
\citep{Kro21} treats density estimation at a fixed point as a special case of a linear functional, and also proposes a Lepski-type procedure for adaptive estimation. Recently, \citep{sart2022density} studied piecewise constant estimators and derived rates of convergence both for the Hellinger and the $L_2$-loss. 

The fundamental problem of universal consistency (both pointwise and with respect to the $L_1$-error), however, has not been considered yet in the context of LDP, and we aim to close this gap in the existing literature. Moreover, we directly derive our results for the multivariate case, which has, with the notable exception of \citep{rohde2020geometrizing}, not been considered yet. In this respect, the derived minimax rates of convergence over Lipschitz classes also complement the existing research literature.

The remaining part of the paper, which is mainly devoted to the proof of analogues of Theorems~\ref{hcons} and \ref{hrate} under LDP, is organized as follows:
first, in Section~\ref{s:hist}, we introduce a privatised version of the multivariate histogram estimator.
Section~\ref{s:results} contains consistency results, namely universal pointwise consistency and strong universal $L_1$-consistency of the private histogram.
In Section~\ref{s:rate} we obtain the rate of convergence of the private histogram for Lipschitz continuous densities, which is complemented by a lower bound in Section~\ref{s:lower}.
Section~\ref{s:simulations} provides a simulation study illustrating our findings and we conclude in Section~\ref{s:discussion} with a brief discussion of our results and comparison with existing approaches.
All proofs are deferred to Appendix~\ref{s:proofs}.

\section{Histogram estimator under LDP}\label{s:hist}

In order to define the privatised histogram estimator we first choose a sequence of closed balls $(S_n)_{n \geq 1} \subset \R^d$ centred at the origin.
The radius of $S_n$ is denoted by $r_n > 0$.
As in the introduction, let $\{A_{h_n,1},A_{h_n,2},\ldots\}$ be a cubic partition of $\R^d$.
Without loss of generality we can assume that the cells $A_{h_n,j}$ are numbered such that $A_{h_n,j}\cap S_n\ne \emptyset$ when $j\le N_n$ for some positive integer $N_n$, and $A_{h_n,j}\cap S_n= \emptyset$ otherwise.
It is not difficult to see that $N_n \asymp (r_n/h_n)^d$.

For nonparametric regression estimation and for classification, \citep{BerrettButucea2019} and \citep{BeGyWa21} introduced a non-interactive privacy mechanism.
In their setup locally privatised samples are created as follows:
the statistician (data base provider) sets and announces some $\sigma_W>0$ and the $i$-th data holder with $i \in \{ 1,\ldots,n \}$ generates and transmits to the statistician the data
\begin{align}\label{eq:def:Z_i}
	Z_{i} \defeq \{W_{n,i,j}, \quad j\le N_n\},
\end{align}
where 
\begin{align}
	\label{Eq:Mech2}
	W_{n,i,j} \defeq \IND_{\{X_i\in A_{h_n,j} \}} + \sigma_W \zeta_{i,j}.
\end{align}
Here the random variables  $\zeta_{i,j}$, $i=1,...,n, \, j=1,2,\dots$ are i.i.d.\,according to a Laplace distribution with unit variance, which has the probability density
\begin{equation*}
	p(x) = \exp(-\sqrt{2}|x|) /\sqrt{2}.
\end{equation*}
Note that  the $i$-th data holder transforms the $d$-dimensional vector $X_i$ into randomised, $N_n$-dimensional data $Z_i$.

Let us briefly recall the definition of $\alpha$-local differential privacy. We refer the reader to the seminal paper by \citep{DuJoWa18} for a comprehensive introduction.
A non-interactive privacy mechanism  is a family of conditional distributions drawing the $i$-th privatised output $Z_i$ from a measurable space $(\mathcal Z,\mathscr Z)$ given the corresponding raw datum $X_i$.
Such a mechanism is said to satisfy the $\alpha$-LDP constraint if
\begin{equation}\label{eq:def:LDP:alpha}
	\sup_{A \in \mathscr Z} \sup_{x,x' \in \mathbb{R}^d} 
	\frac{\PROB(Z_i\in A \mid X_i = x)}{\PROB(Z_i\in A \mid X_i = x')} 
	\leq e^\alpha, \quad i=1,2,\dots   
\end{equation}
Here, $\alpha$ is a non-negative privacy parameter with the interpretation that smaller values of $\alpha$ lead to a stronger privacy guarantee.

Standard calculations show that the mechanism defined via \eqref{eq:def:Z_i} and \eqref{Eq:Mech2} satisfies the LDP constraint provided that the standard deviation $\sigma_W = \sigma_W(\alpha)$ is chosen sufficiently large, namely $\sigma_W \geq 2^{3/2}\alpha^{-1}$, see \citep{BerrettButucea2019} and  \citep{BeGyWa21}.
In the sequel, we assume equality such that
\begin{align}
	\label{sa}
	\sigma_W = 2^{3/2}\alpha^{-1}
\end{align}
holds.

Though the privacy mechanism defined through \eqref{eq:def:Z_i} and \eqref{Eq:Mech2} has been suggested in others papers already for univariate density estimation, we propose a novel estimator based on the multivariate privatised data $Z_i$, $i=1,\ldots,n$.

\citep{DuJoWa18} proposed the estimator $\hat f_n$ defined by
\begin{align*}
	\hat f_{n}(\xb)= \frac{\hat\mu_n(A_{h_n,j})}{h_n^d} \,\IND_{\{j\le N_n \}}\qquad \mbox{if } \xb\in A_{h_n,j},
\end{align*}
where
\begin{equation}\label{eq:alternative:muhat}
	\hat \mu_n(A_{h_n,j}) 
	= \frac{1}{n} \sum_{i=1}^n W_{n,i,j}
	=\nu_n(A_{h_n,j})+\mu_n(A_{h_n,j})
\end{equation}
with
\begin{align*}
	\nu_n(A_{h_n,j})=\frac{\sigma_W}{n} \sum_{i=1}^n \zeta_{i,j}.
\end{align*}
Note that the $\hat \mu_n(A_{h_n,j})$ are unbounded and not necessarily non-negative.

We now introduce a novel estimate, which improves upon $\hat f_n$.
Using the privatised data in~\eqref{Eq:Mech2} one can compute the empirical distribution function
\begin{equation*}
G_{n,j}(z)\defeq \frac 1n \sum_{i=1}^n\IND_{\{W_{n,i,j}\le z\}}
\end{equation*}
for all $j\leq N_n$.
One has
\begin{align*}
	G_{n,j}(z)
	&=\frac 1n  \sum_{i=1}^n\IND_{\{\IND_{\{X_i\in A_{h_n,j} \}} + \sigma_W \zeta_{i,j}\le z\}}\\
	&=\frac 1n  \sum_{i=1}^n(\IND_{\{X_i\in A_{h_n,j} \}}\IND_{\{1 + \sigma_W \zeta_{i,j}\le z\}}
	+(1-\IND_{\{X_i\in A_{h_n,j} \}})\IND_{\{\sigma_W \zeta_{i,j}\le z\}})\\
	&=\frac 1n  \sum_{i=1}^n(\IND_{\{X_i\in A_{h_n,j} \}}(\IND_{\{1 + \sigma_W \zeta_{i,j}\le z\}}-\IND_{\{\sigma_W \zeta_{i,j}\le z\}})
	+\IND_{\{\sigma_W \zeta_{i,j}\le z\}}).
\end{align*}
Denoting with $H$ the standard Laplace distribution function, conditioning on the raw data $\mathcal D_n = \{ X_1,\ldots,X_n \}$ yields
\begin{align*}
	&\EXP[G_{n,j}(z)\mid \mathcal D_n]\\
	&=
	\frac 1n  \sum_{i=1}^n (\IND_{\{X_i\in A_{h_n,j} \}}(H((z-1)/\sigma_W)-H(z/\sigma_W))+H(z/\sigma_W))\\
	&=
	H(z/\sigma_W) + \mu_n(A_{h_n,j})(H((z-1)/\sigma_W)-H(z/\sigma_W)),
\end{align*}
which can be rearranged to obtain
\begin{align}
	\label{Eq:MG}
	\mu_n(A_{h_n,j})
	&=\frac{H(z/\sigma_W)-\EXP\{G_{n,j}(z)\mid \mathcal D_n \}}{H(z/\sigma_W)-H((z-1)/\sigma_W)}.
\end{align}
Based on \eqref{Eq:MG}, we set
\begin{align*}
	\tilde \mu_n(A_{h_n,j})
	&=\frac{H(z/\sigma_W)-G_{n,j}(z)}{H(z/\sigma_W)-H((z-1)/\sigma_W)}
\end{align*}
as an estimator of $\mu_n(A_{h_n,j})$, and finally define the corresponding density estimator from the privatised data $Z_1,\ldots,Z_n$
\begin{align*}
\ftilde_{n}(\xb)= \frac{\tilde\mu_n(A_{h_n,j})}{h_n^d} \,\IND_{\{j\le N_n \}}\qquad \mbox{if } \xb\in A_{h_n,j}.
\end{align*}
Then,
\begin{align*}
	\Var (\tilde \mu_n(A_{h_n,j}))
	=\frac{\Var (G_{n,j}(z))}{(H(z/\sigma_W)-H((z-1)/\sigma_W))^2},
\end{align*}
and the denominator on the right-hand side is maximized by taking $z=1/2$ which we consider from now on only.
For this choice of $z$, $\tilde \mu_n(A_{h_n,j})$ can be expressed as
\begin{align}
	\label{Eq:TMG}
	\tilde \mu_n(A_{h_n,j})
	&=\frac{H(1/(2\sigma_W))-G_{n,j}(1/2)}{H(1/(2\sigma_W))-H(-1/(2\sigma_W))}\nonumber\\
	&=\frac{H(1/(2\sigma_W))-G_{n,j}(1/2)}{1-2H(-1/(2\sigma_W))}.
\end{align}

In contrast to the classical histogram estimator $f_n$ the privatised histogram estimator $\tilde f_n$ is not necessarily a probability density function: neither positivity, nor the property that the estimator integrates to one are in general satisfied.
For this reason we may modify this estimator by taking its positive part and then normalise such that the modified estimator integrates to $1$:
\begin{equation}
	\label{eq:def:fhat:mod}
	\ftilde_{n,\mathrm{mod}}(\xb) = \frac{\ftilde_n(\xb) \vee 0}{\int_{\mathbb R^d}  (\ftilde_n(\zb) \vee 0) \dd \zb}.
\end{equation}
$\ftilde_{n,\mathrm{mod}}$ is called the positive projection of $\ftilde_{n}$.
It is well-known (see \citep{DeGy85}, pp.~269--70) that the positive projection always reduces the $L_1$-error of the estimator:
\begin{align*}
	\int_{\R^d} |f(\xb) - \ftilde_{n,\mathrm{mod}}(\xb) | \dd \xb
	\le
	\int_{\R^d} |f(\xb) - \ftilde_n(\xb) | \dd \xb.
\end{align*}
Because of
\begin{align*}
	\EXP[\tilde \mu_n(A_{h_n,j})\mid \D_n]
	&=
	\EXP[\hat \mu_n(A_{h_n,j})\mid \D_n]
	=\mu_n(A_{h_n,j}),
\end{align*}
the bias of both density estimators $\tilde f_n$ and $\hat f_n$ coincides with the one of the non-private histogram $h_n$, and the results obtained in the sequel for the estimator $\ftilde_n$ hold equally true for the alternative estimator $\hat f_n$ with only minor modifications necessary in the proofs but the same requirements concerning the choice of the bandwidth parameter $h_n$.
However, the less intuitive definition of $\tilde \mu_n(A_{h_n,j})$ in \eqref{Eq:TMG} outperforms the obvious estimator $\hat \mu_n(A_{h_n,j})$ in \eqref{eq:alternative:muhat} both in finite sample studies and also from a theoretical point of view (see Remark~\ref{remark} below).
Numerical results confirming this will be presented in Section~\ref{s:simulations}.

\section{Universal consistency}\label{s:results}

Our first theoretical result establishes the strong universal $L_1$- and pointwise  consistency of the privatised histogram and can be seen as an LDP-analogue of Theorem~\ref{hcons}.
\begin{theorem}
	\label{Thm:pL1}
	Assume that $r_n \to \infty$.
	If $h_n \to 0$ and $nh_n^{2d}/\log n \to \infty$, then
	\begin{equation*}
		\lim_{n \to \infty} \ftilde_{n} (\xb) 
		= f(\xb) \quad \text{almost surely}
	\end{equation*}
	for Lebesgue-almost every $\xb \in \R^d$. 
	If, in addition, $nh_n^{2d}/(r_n^{2d} \log n) \to \infty$, then
	\begin{equation*}
\lim_{n\to \infty} \int_{\R^d} | \ftilde_n(\xb) - f(\xb) | \dd \xb
		= 0 \quad \text{almost surely.}
	\end{equation*}
\end{theorem}

First note that a similar theorem can be proven for $\hat f_n$ instead of $\tilde f_n$.
We remark that the proof of the pointwise consistency in Section~\ref{subs:proof:PUC} can equally be used to prove universal pointwise consistency of the non-private histogram $f_n$.
In that case one can use Bennett's inequality (see \citep{boucheron2013concentration}, Theorem~2.9) instead of Hoeffding's inequality, and the (weaker) conditions to obtain universal pointwise consistency in that case turn out to be $h_n \to 0$ and $nh_n^{d}/\log n \to \infty$. 
This result seems to be novel as well and might be of independent interest, since establishing universal pointwise consistency usually comes along with quite technical assumptions, see \citep{vidal2005pointwise}.

A direct application of Theorem~\ref{Thm:pL1} can be derived from the fact that $L_1$-consistent density estimates yield distribution estimates that are consistent in total variation. 
More precisely, let $\mu_n^*$ be the distribution estimate derived from the histogram estimate $\tilde f_n$ as follows:
\begin{align*}
	\mu_n^*(A)
	&=\int_A \ftilde_n(\xb) \dd \xb
	=\sum_{j=1}^{N_n}\int_{A\cap A_{h_n,j}} \ftilde_n(\xb) \dd \xb
	=\sum_{j=1}^{N_n}\frac{\lambda(A\cap A_{h_n,j})}{\lambda(A_{h_n,j})} \tilde\mu_n(A_{h_n,j}),
\end{align*}
where $\lambda$ is the Lebesgue measure.
Thus, under the conditions of Theorem \ref{Thm:pL1} one gets
\begin{align*}
	\TV(\mu,\mu_n^*) \defeq \sup_{A\in \mathscr B}|\mu(A)-\mu^*_n(A)| = \frac 12 \int_{\R^d} |\ftilde_n(\xb) - f(\xb)| \dd \xb
	\to 0
\end{align*}
almost surely, where $\mathscr B$ stands for the Borel $\sigma$-algebra.
If the distribution $\mu$ is a mixture of absolutely continuous and discrete distributions, then from non-private data, \citep{BaGyvdM92} introduced a distribution estimate that is consistent in total variation.
It remains an open question whether consistent estimation in total variation distance is possible from private data without knowing the support of the discrete component.

\section{Rate of convergence}\label{s:rate}

Our next theorem provides an LDP analogue of Theorem~\ref{hrate} and states a weak consistency rate for Lipschitz functions with compact support.
\begin{theorem}
	\label{Thm:prate}
	Under the conditions of Theorem \ref{hrate}, one has that
	\begin{align*}
		&\EXP\left[  \int_{\R^d} |f(\xb ) - \ftilde_n(\xb) | \dd \xb\right]\\
		&\leq
		L\sqrt{d}h_n + \int_{\R^d \setminus \bigcup_{j=1}^{N_n} A_{h_n,j}} f(\xb) \dd \xb
		+ O \left( \frac{r_n^d}{\sqrt{nh_n^d}} \right) + O \left( \frac{\sigma_Wr_n^d}{\sqrt{nh_n^{2d}}} \right).
	\end{align*}
	Assume that $r_n \to \infty$. Then, for $n$ sufficiently large, $S \subseteq S_n$, and the bound becomes
	\begin{align}\label{eq:bound:upper:O}
		\EXP\left[ \int_{\R^d} |f(\xb) - \ftilde_n(\xb) | \dd \xb\right] 
		&\leq
		L\sqrt{d}h_n + O \left( \frac{r_n^d}{\sqrt{nh_n^d}} \right) + O \left( \frac{\sigma_Wr_n^d}{\sqrt{nh_n^{2d}}} \right).
	\end{align}
	Again, the same guarantees hold for $\tilde f_n$ being replaced with $\hat f_n$.
\end{theorem}

\begin{remark}\label{remark}
	The proof of Theorem~\ref{Thm:prate} is based on bounding the $L_1$-error by the sum of an approximation error (bias) and a stochastic error term.
	As mentioned already above the bias of both private estimators $\ftilde_n$ and $\hat f_n$ is equal to the one of the classical histogram $f_n$.
	However, the proof of Theorem~\ref{Thm:prate} shows that
	\begin{equation*}
		\Var (\tilde \mu_n(A_{h_n,j}))
		\le 
		\Var (\hat \mu_n(A_{h_n,j}))
	\end{equation*}
	yielding a bound for the stochastic error for $\ftilde_n$ that is less than or equal to the one obtained for $\hat f_n$.
	More precisely, from the proof of Theorem~\ref{Thm:prate} one can even derive an upper bound of order $e^{-1/(\sqrt 2 \sigma_W)}/n$ (when $\sigma_W \to 0$) for the quantity $\Var(\tilde \mu_n(A_{h_n,j}))$ which is faster than the rate $\sigma_W^2/n$ obtained for the quantity $\Var(\hat \mu_n(A_{h_n,j}))$.
	Note that the minimax lower bound suggests a lower bound of order $e^{-2\alpha}/n$.
	Hence, our modified estimator provides an essential improvement in the regime where $\sigma_W$ is relatively small towards the standard estimator $\hat f_n$.
	We emphasize that one also has to consider a regime where the non-private rate is dominating when $\sigma_W \to 0$ resp. $\alpha \to \infty$ as $n \to \infty$.
	The exact expression of the optimal rate of convergence that is sharp with respect to both the privacy level $\alpha$ and numerical constants in all possible regimes, however, remains an important open theoretical question which is outside the scope of the present paper. 
\end{remark}

In practice, one will usually consider a fixed privacy level $\alpha$.
Consequently, $\sigma_W$ is also fixed and the third term on the right-hand side of \eqref{eq:bound:upper:O} dominates the second one for $n$ being sufficiently large.
In this case the optimal choice of the bandwidth is $h_n=c\cdot (n/(\sigma_W^2r_n^{2d}))^{-1/(2d+2)}$, and \eqref{eq:bound:upper:O} yields
\begin{align}
	\label{as}
	\EXP\left[ \int_{\R^d} |f(\xb) - \tilde f_n(\xb) | \dd \xb\right] 
	\le
	(n\alpha^2/r_n^{2d})^{-1/(2d+2)}.
\end{align}

The quantity $r_n$ in the rate of convergence of Theorem~\ref{Thm:prate} is caused by the fact that all cells $A_{h_n,j}$ with $j=1,\ldots,N_n$ contribute to the private histogram estimator regardless of whether the cell actually has non-zero mass under the distribution $\mu$.
On the contrary, when the raw data $\mathcal D_n$ are accessible, the classical histogram can, with probability one, only put mass on cells that have non-empty intersection with the support of the underlying distribution.
If a priori knowledge on the support $S$ of the density was given in the LDP framework (which is rather unrealistic), one could also consider only those cells that have non-empty intersection with $S$, and the contribution of $r_n$ could be avoided.

We remark that Theorem~\ref{Thm:prate} can easily be extended from Lipschitz to more general H\"older continuous functions: in the proof of this extended result, only the bound for the bias term would change.
Let us mention that such general rates of convergence for nonparametric density estimation under LDP have already been derived in the literature in several papers \citep{DuJoWa18,BuDuKrSa20,sart2022density}.
However, these articles are restricted to the case of univariate densities only, and our results can be seen as a first step towards an extension to the multivariate setup.
The rate derived in Theorem~\ref{Thm:prate} indicates a worsening of the curse of dimensionality under LDP by a factor of $2$ in the exponent.
This phenomenon has already been established by \citep{rohde2020geometrizing} for the estimation of linear functionals under local differential privacy.

\section{Minimax lower bound}\label{s:lower}

The purpose of this section is to prove that the rate obtained in Theorem~\ref{Thm:prate} is essentially optimal, that is, no other privacy mechanism together with an accompanying estimator respecting the $\alpha$-LDP constraint can attain a faster rate of convergence.
Here, in accordance with the approach in \citep{DuJoWa18}, we will even allow for a richer class of admissible privacy mechanisms than only the non-interactive ones introduced in Section~\ref{s:hist}.
The amount of interaction that is generally allowed under LDP assumes that the individual data holders can generate their private data in dependence on both their own sensitive raw data and the anonymised data of other data holders that have been generated before. In mathematical terms, this corresponds to the existence of Markov kernels
\begin{equation}\label{eq:Q:markov:kernel}
	Q_i\colon \Zs_i \times (\R^d \times \Zc_1 \times \ldots \times \Zc_{i-1}) \to [0,1]
\end{equation}
for $i=1,\ldots,n$ where $(\Zc_i,\Zs_i)$ denotes the measure space for the private data of the $i$-th data holder.
The amount of interaction allowed under LDP makes a distinction between LDP and the notion of global differential privacy (GDP). Under GDP, the data holders send their raw data to a trusted central curator. This trusted curator can perform arbitrary analyses with the raw data, and only the output published by this trusted entity has to be anonymised. This notion is less restrictive then the notion of LDP since the trusted curator can reproduce any procedure that is admissible under LDP.

With the general definition of admissible privacy mechanisms under LDP, condition \eqref{eq:def:LDP:alpha} has to be replaced with
\begin{equation}\label{eq:def:DP:alpha:general}
	\sup_{A \in \mathscr Z} \sup_{\substack{x,x' \in \mathbb{R}^d\\z_j \in \Zc_j}} 
	\frac{Q_i(A \mid X_i = x,Z_1=z_1,\ldots,Z_{i-1}=z_{i-1})}{Q_i (A \mid X_i = x',Z_1=z_1,\ldots,Z_{i-1}=z_{i-1})} 
	\leq e^\alpha
\end{equation}
for $i=1,\ldots,n$.
Let us denote with $\Qc_\alpha$ the set of all privacy mechanisms $Q=(Q_1,\ldots,Q_n)$ with $Q_i$ as in \eqref{eq:Q:markov:kernel} and satisfying \eqref{eq:def:DP:alpha:general}.
The quantity of interest is the private minimax risk over a class $\Fc$ of densities which is given by
\begin{equation*}
	\inf_{\substack{\ftilde \\ Q \in \Qc_\alpha}} \sup_{f \in \Fc} \EXP_f \left[ \int_{\R^d} \lvert \ftilde(\xb) - f(\xb) \vert \dd\xb \right],
\end{equation*}
where the infimum is taken over all privacy mechanisms and estimators $\ftilde$ defined in terms of the anonymised data $Z=(Z_1,\ldots,Z_n)$.

For a constant $L > 0$, we now consider the \emph{Lipschitz class} $\Flipd(L)$ defined as the class of all densities $f\colon \R^d \to \R$ vanishing outside $[0,1]^d$ and satisfying the Lipschitz condition
\begin{equation*}
	\rvert f(\xb) - f(\yb) \lvert \leq L \lVert \xb - \yb \rVert,
\end{equation*}
where $\lVert \, \cdot \, \rVert$ denotes the Euclidean norm of $\R^d$.

The following theorem complements the convergence rate derived in Theorem~\ref{Thm:prate} by a lower bound.
Its proof essentially relies on an appropriate information theoretical inequality between Kullback-Leibler divergence of the privatised data on the one side and the squared total-variation distance of the raw data on the other side, which has been derived in \citep{DuJoWa18} and is stated in Section~\ref{s:proofs}.

\begin{theorem}[Private lower bound for $\Flipd$]\label{THM:LOWER:LIPSCHITZ} For all $L$ sufficiently large, we have
	\begin{equation*}
		\inf_{\substack{\ftilde \\ Q \in \Qc_\alpha}} \sup_{f \in \Flipd(L)} \EXP_f \left[ \int_{\R^d} \lvert \ftilde(\xb) - f(\xb) \vert \dd\xb \right] \gtrsim (n(e^{\alpha}-1)^2)^{-\frac{1}{2d + 2}},
	\end{equation*}
	where the infimum is taken over all, potentially interactive, privacy channels $Q$ ensuring $\alpha$-LDP, and all estimators $\ftilde$ based on the privatised observations $Z_1,\ldots,Z_n$.
\end{theorem}

As for Theorem~\ref{Thm:prate}, one can rewrite this result for the interesting privacy regime where $\alpha\in [0,1]$. Then the lower bound in Theorem \ref{THM:LOWER:LIPSCHITZ} can be rewritten as
\begin{align}
	\label{ass}
	(n\alpha^2)^{-1/(2d+2)}.
\end{align}
Thus, the upper bound (\ref{as}) and the lower bound (\ref{ass}) are matching, 
i.e., the obtained rate (\ref{as}) is essentially optimal and cannot be improved by any, potentially interactive, privacy mechanism.
We emphasise that the privacy mechanism defined by~\eqref{eq:def:Z_i} and \eqref{Eq:Mech2} is of the preferential, non-interactive, form.

Using essentially the same techniques as the ones used in the proof of Theorem~\ref{THM:LOWER:LIPSCHITZ}, private minimax lower bounds for more general smoothness classes as considered in \citep{DeGy85} can be derived.
More precisely, for Hölder smooth densities with index $\beta\in (0,1]$ one can establish the lower bound $n^{-\beta/(2\beta+2d)}$.

\section{Simulation study}\label{s:simulations}

\begin{figure}
	\centering
	\begin{subfigure}[b]{0.45\textwidth}
		\centering
		\includegraphics[width=\textwidth]{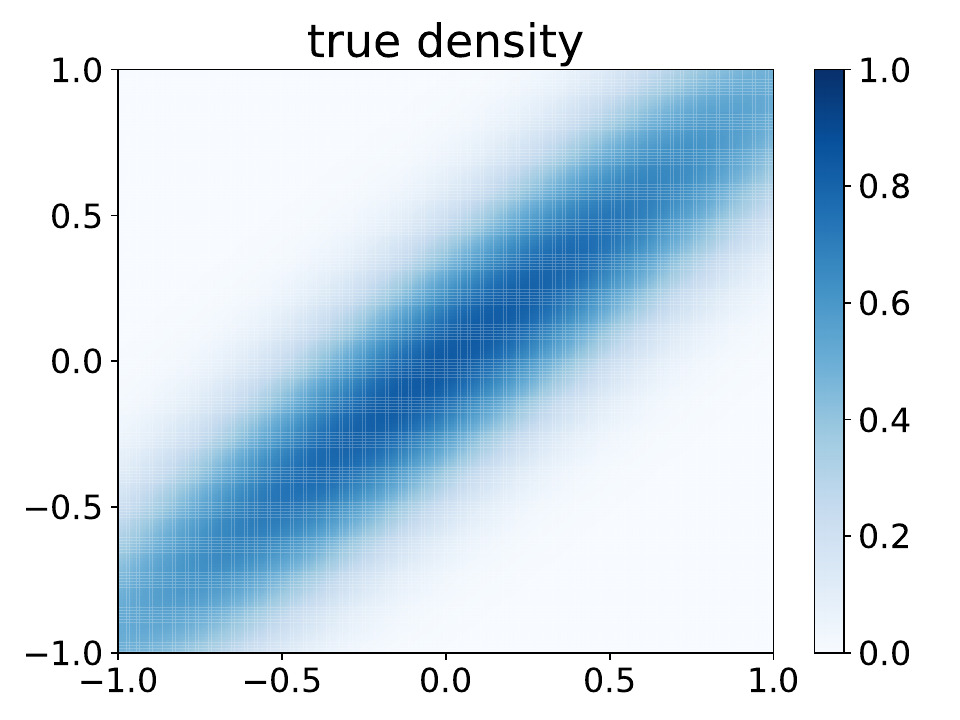}
		\caption{}
\end{subfigure}
	\hfill
	\begin{subfigure}[b]{0.45\textwidth}
		\centering
		\includegraphics[width=\textwidth]{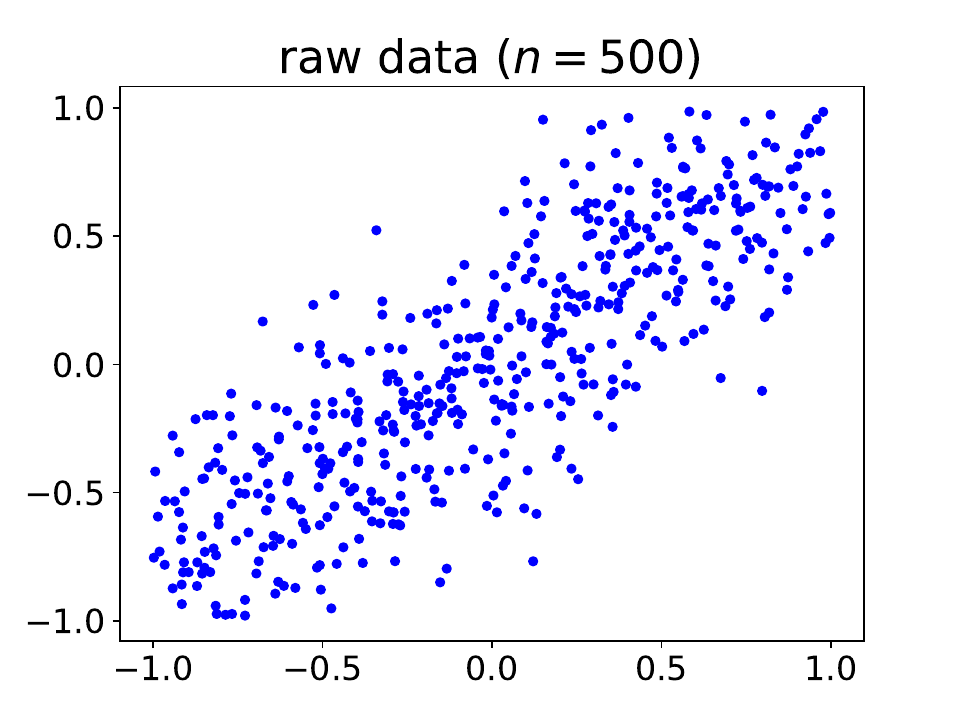}
		\caption{}
\end{subfigure}
	\caption{(a) True probability density function corresponding to a correlated truncated two-dimensional Gaussian distribution. (b) i.i.d.\,sample of size $n=500$ from the density function in~(a).}
	\label{fig:true_density_and_sample}
\end{figure}

We compared the different versions of the privatised histogram estimator in a small simulation study.
Results of these experiments are reported in this section.
We considered the two-dimensional case where the raw data were generated from a truncated two-dimensional Gaussian.
More precisely, $X_1,\ldots,X_n$ were drawn i.i.d.\,$\sim \mathcal N(\boldsymbol 0, \Sigma)$ restricted to the square $[-1,1]^2$ where the covariance matrix was chosen as
\begin{equation*}
	\Sigma = \begin{pmatrix}
		1 & 0.9\\
		0.9 & 0.9
	\end{pmatrix}.
\end{equation*}
The true density function of this data and a typical raw data sample of size $n=500$ are plotted in Figure~\ref{fig:true_density_and_sample}.

In light of the drastic sample size reduction under LDP, we used a large sample size of $n=10^5$ for the simulations.
We considered different values for the privacy parameter $\alpha$, namely $\alpha \in \{ \infty,0.5,0.25,0.1 \}$.
Here, setting $\alpha = \infty$ corresponds to the case without privacy constraints, and in this case we consider the classical histogram estimator $f_n$ defined in \eqref{hist}.
Recall that by Equation~\eqref{sa} the privacy parameter $\alpha$ is related to the variance $\sigma^2_W$ appearing in the perturbation of the histogram via the identity $\sigma_W^2 = 8/\alpha^2$.
As a consequence, note that even in the case of the weakest privacy guarantee with $\alpha = 0.5$, a considerable amount of Laplace noise with standard deviation $\sigma_W \approx 5.66$ has to be added which explains the difficulty of inference under LDP.

In order to further illustrate this effect we compare the information stored in the classical histogram of any data holder with the one of the corresponding perturbed histogram.
The non-privatised histogram generated locally by any data holder can be stored as a binary matrix with exactly one entry equal to $1$ indicating the cell, where the respective data holder's observation is located.
An example of such a matrix is given by
\begin{equation*}\footnotesize
	\begin{pmatrix}
		0 & 0 & 0 & 0 & 0\\
		0 & 1 & 0 & 0 & 0\\
		0 & 0 & 0 & 0 & 0\\
		0 & 0 & 0 & 0 & 0\\
		0 & 0 & 0 & 0 & 0
	\end{pmatrix}.
\end{equation*}
The strong effect of the anonymization procedure becomes apparent, when considering the perturbed version of this matrix, where we have taken $\sigma_W^2 = 32$.
A realisation of such a perturbed matrix containing the $W_{n,i,j}$ in \eqref{Eq:Mech2} for a fixed $i$ is given by
\begin{equation*}\footnotesize
	\begin{pmatrix}
		-0.18 & -1.99 &  2.64 &  0.83 & -8.85\\
		4.36 & 0.58 & 2.59  & -1.66 & 2.17\\
		2.01 & -2.49 & -0.83 &  8.20 & -0.50\\
		0.12 & -0.84 & -13.20 & 2.15 & 6.01\\
		-3.97 & 4.39 & 2.45 & 4.28 & 1.51
	\end{pmatrix}.
\end{equation*}

\begin{figure}[htb]
	\centering
	\includegraphics[width=1.0\textwidth]{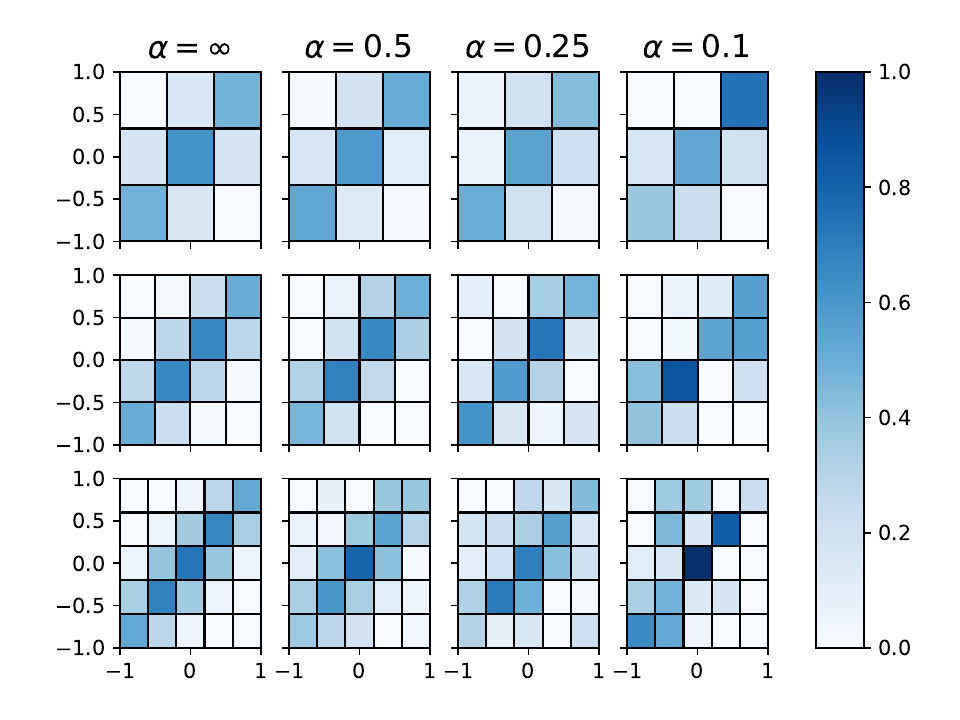}
	\caption{Illustration of the projected histogram estimator $\tilde f_{n, \mathrm{mod}}$ of the true density $f$ shown in Figure~\ref{fig:true_density_and_sample} for different cubic partitions of the domain $[-1,1]^2$, and different privacy levels $\alpha$. The level of privacy imposed via LDP increases from left to right whereas accuracy of the estimates decreases at the same time.
	}\label{fig:3-5}
\end{figure}

For our simulations, we split the domain $[-1,1]^2$ into $k \times k$-grids of squares for $k = 3,4,5$. This corresponds to choices of the bandwidth parameter $h_n$ equal to $2/3$, $1/2$, and $2/5$, respectively.
The resulting estimators are shown in Figure~\ref{fig:3-5}. In all of these figures the privatised histogram $\tilde f_{n,\mathrm{mod}}$ is plotted on the square $[-1,1]^2$.
Already this limited number of simulations confirms the general message from the theoretical results that a smaller value of $\alpha$, which yields more privacy, must be paid for by a coarser resolution, i.e., a larger bandwidth in contrast to the non-privacy framework.
Whereas for the $3 \times 3$-grids the trend of the data concentrating close to the diagonal can be detected for all values of $\alpha$, this trend becomes more vague for the $4 \times 4$ and $5 \times 5$ grid.
For the $5 \times 5$ grid the choice $\alpha = 0.5$ is still close to the non-privatized histogram, and the estimator for $\alpha = 0.1$ now completely fails to detect the trend in the data.

\FloatBarrier

\begin{figure}[htb]
	\centering
	\includegraphics[width=0.85\textwidth]{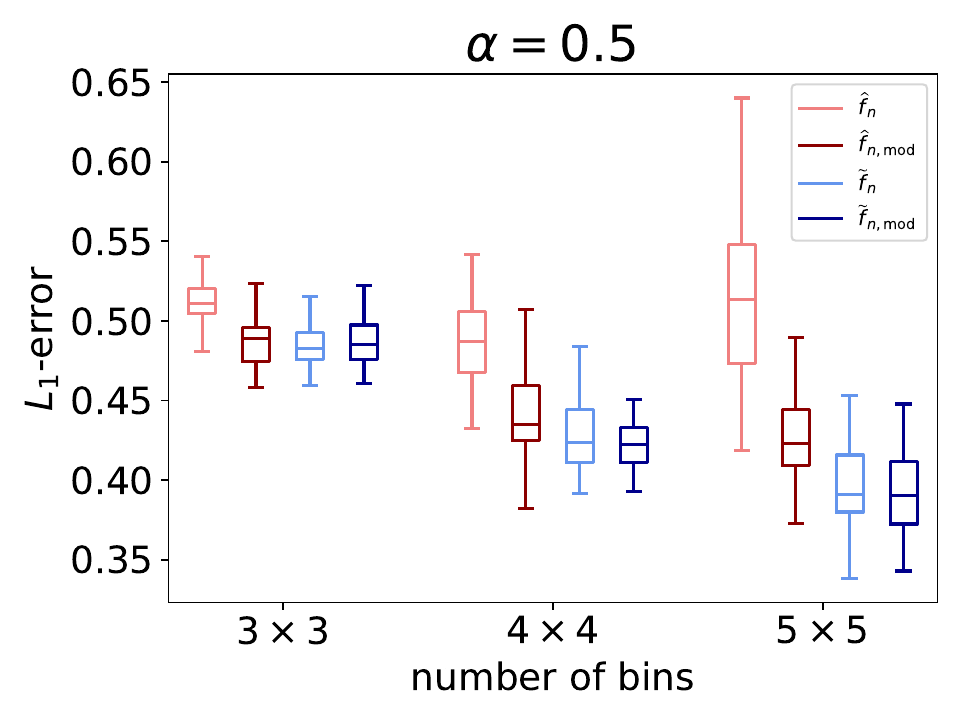}
	\caption{$L_1$-error in dependence on the chosen partition and method for $\alpha = 0.5$. The boxplots are based on a Monte-Carlo study with $50$ iterations with $n=10^5$ for each setup.}
	\label{fig:L1_error_alpha_0_5}
\end{figure}

\begin{figure}[htb]
	\centering
	\includegraphics[width=0.85\textwidth]{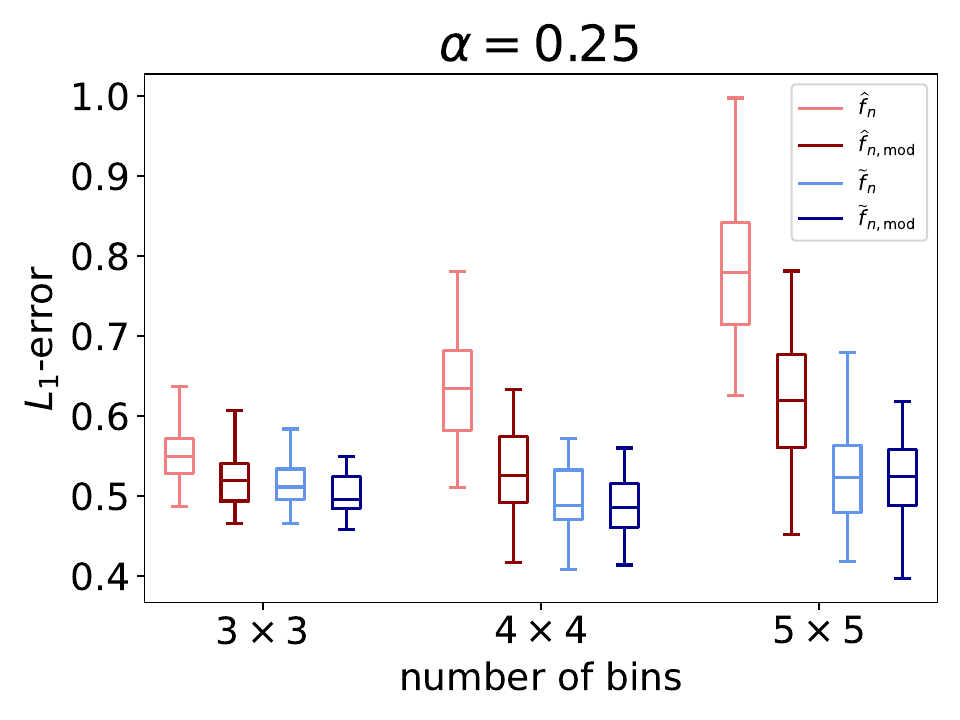}
	\caption{$L_1$-error in dependence on the chosen partition and method for $\alpha = 0.25$. The boxplots are based on a Monte-Carlo study with $50$ iterations with $n=10^5$ for each setup.}
	\label{fig:L1_error_alpha_0_25}
\end{figure}

Figures~\ref{fig:L1_error_alpha_0_5} and \ref{fig:L1_error_alpha_0_25} compare the $L_1$-error for the standard histogram $f_n$ and our novel estimator $\tilde f_n$ (both with and without the projection given through \eqref{eq:def:fhat:mod}) for $\alpha = 0.5$ and $\alpha = 0.25$ and different partitions.
The figures indicate that the novel estimator outperforms the obvious one and also show the positive effect of the projection \eqref{eq:def:fhat:mod}.
Interestingly, we detected a monotone improvement in the order $\int |\hat f_n -f|$, $\int |\hat f_{n, \mathrm{mod}} -f|$, $\int |\ftilde_n -f|$, $\int |\ftilde_{n, \mathrm{mod}} -f|$.
Furthermore, the projection is more effective for $\hat f_{n, \mathrm{mod}}$, because $\hat f_{n}$ takes negative values with large probability.
In addition, the figures demonstrate the huge impact of LDP (note that different scalings of the $y$-axis are used for the two figures).

\FloatBarrier

\section{Discussion}\label{s:discussion}

In this paper, we have proven strong universal consistency, both pointwise and in the $L_1$-sense, for an appropriately specified histogram density estimator under LDP for multivariate raw data.
We have also derived minimax rates of convergence for Lipschitz continuous densities.
It is certainly of interest to compare the obtained results under LDP with the well-known results in the classical case when no privacy restrictions are imposed.
Such a comparison is given in Table~\ref{tab:summary} in condensed form.
For instance, in order to obtain universal pointwise consistency the prerequisites for the private histogram ($h_n \to 0, r_n \to \infty, nh_n^{2d}/\log n \to \infty$) are harder than those for the private one ($h_n \to 0, nh_n^d/\log n \to \infty$), and similarly for strong universal consistency (see row (SUC) in Table~\ref{tab:summary}).
Moreover, the rate of convergence for Lipschitz continuous functions degrades from $n^{-1/(d+2)}$ in the classical case to $n^{-1/(2d+2)}$ in the LDP case.

\renewcommand{\arraystretch}{1.45}
\begin{table}
	\centering
	\begin{tabular}{c!{\color{black}\vrule}c!{\color{gray!50}\vrule}c}
		& \footnotesize non-private histogram &\footnotesize privatised histogram\\
		\arrayrulecolor{black}\hline
		\scriptsize (UPC) & $h_n \to 0$, $nh_n^d/\log n \to \infty$ & $h_n \to 0$, $r_n \to \infty$, $nh_n^{2d}/\log n \to \infty$\\
		\arrayrulecolor{gray!20}\hline
		\scriptsize (SUC) & $h_n \to 0, nh_n^d \to \infty$ & (UPC) + $nh_n^{2d}/(r_n^{2d} \log n) \to \infty$\\
		\arrayrulecolor{gray!20}\hline
		\scriptsize (LIP) & $n^{-\frac{1}{d+2}}$ &  $n^{-\frac{1}{2d+2}}$ \\
	\end{tabular}
	\caption{Comparison of classical and (modified) privatised histogram estimator.
		The rows indicate: (UPC) assumptions for universal pointwise consistency, (SUC) assumptions for strong universal consistency, and (LIP) minimax optimal rate over Lipschitz classes.}
	\label{tab:summary}
\end{table}
In this work, we restricted ourselves to histogram estimators of the density.
This choice was mainly motivated by our primary interest in establishing strong universal consistency where results for histogram estimators are particularly easy to derive.
Deriving universal consistency under LDP for other estimators like wavelet estimators (see \citep{chacon2005consistency} for the classical case) or general delta estimators as in \citep{vidal2005pointwise} is an interesting direction for future research.
Note that kernel density estimators, although well-understood in the classical case (see \citep{DeGy85}), are more difficult to treat in the context of LDP.
Note that from our privatised data one cannot compute a standard kernel type density estimator.
Whereas anonymization in the LDP sense of kernel estimators at one fixed point is possible by means of Laplace perturbation (see \citep{Kro21}), recent results in \citep{reimherr2019elliptical} suggest that anonymization of the whole kernel density estimator cannot be achieved by a mere perturbation approach, at least within a natural class of potential perturbation distributions.

Finally, we mention that our results on minimax convergence rates can be extended to include Hölder smoothness with index $\beta \leq 1$ as well, and also anisotropy can be included.
In this regard, the question of adaptation also arises naturally, and this issue certainly deserves to be considered in more detail in future work.
It also remains an open problem whether data generated by our privacy mechanism can be used for an estimator having faster convergence rate when the true density has smoothness $\beta>1$.

Our final open question concerns the properties of the centred $L_1$-error  $U_n=\int |\ftilde_n -f|-\EXP[\int |\ftilde_n -f|]$ under LDP.
Deriving the rate at which the centred error $U_n$ converges to zero and asymptotic normality of an appropriately scaled version of $U_n$ provide interesting problems for future work.
This latter problem could be seen as an extension of the results derived in \citep{berlinet1995asymptotic} to the LDP setup.

\section*{Acknowledgements}

We thank two anonymous referees for their careful reading and suggestions that led to a substantial improvement of the paper.

\section*{Funding information}

The research of Martin Kroll was supported by the German Research Foundation (DFG) under the grant DFG DE 502/27-1.
The research of Martin Kroll was done at the Ruhr-Universität Bochum, final revisions were done at Universität Bayreuth.

\printbibliography

\newpage

\appendix

\section{Proofs}\label{s:proofs}
\subsection{Proof of Theorem~\ref{Thm:pL1}}
\label{subs:proof:PUC}

The following proof of universal consistency properties of the private histogram estimators relies on concentration inequalities for quantities that depend on the perturbation terms $W_{n,i,j}$. The necessary amount of concentration leads to conditions for universal consistency that are stronger than for classical histograms defined in terms of the raw data $\mathcal D_n$.
First, we prove the pointwise consistency.
The proof is based on the decomposition
\begin{equation}\label{eq:bias:var:pointwise}
	\lvert \ftilde_n(\xb) - f(\xb) \rvert \leq \lvert \EXP\ftilde_n(\xb) - f(\xb) \rvert + \lvert \ftilde_n(\xb) - \EXP \ftilde_n(\xb) \rvert.
\end{equation}
First, for all $n \geq n_0$ with $n_0=n_0(\xb)$ chosen sufficiently large we have $\xb \in S_n$.
For such $n$, the identity
\begin{equation*}
	\EXP \ftilde_n(\xb) = \EXP f_n(\xb) = \frac{\mu(A_{n}(\xb))}{h_n^d}
\end{equation*}
holds, where we denote with $A_n(\xb)$ the unique cube from the collection of $A_{h_n,j}$ that contains the considered $\xb$.
Then, the consistency of the bias term (the first term on the right-hand side of \eqref{eq:bias:var:pointwise}) follows from the generalised Lebesgue density theorem (cf. \cite{wheeden2015measure}, Theorem~7.16) as follows: let $B_n(\xb)$ denote the smallest cube centred at $\xb$ containing $A_n(\xb)$.
Again, let $\lambda$ denote the Lebesgue measure.
If there is a universal constant $c>0$ such that
\begin{equation}\label{eq:cube:condition}
	\lambda(B_n(\xb)) \leq c\lambda(A_n(\xb)),
\end{equation}
then $h_n \to 0$ implies $\lvert \EXP\ftilde_n(\xb) - f(\xb) \rvert \to 0$ for $\lambda$-almost every $\xb$.
Obviously, \eqref{eq:cube:condition} is satisfied with $c=2^d$.

In order to bound the second term on the right-hand side of \eqref{eq:bias:var:pointwise}, we use Hoeffding's inequality (in the formulation of \citep{boucheron2013concentration}, Theorem~2.8) which can be applied, since the random variables $W_{n,i,j}$ for $i=1,\ldots,n$ and \emph{fixed $j$} are independent.
Denoting with $j=j(n)$ the index such that $A_n(\xb) = A_{h_n,j}$ we can thus conclude
\begin{align*}
	&\PROB(\lvert \ftilde_n(\xb) - \EXP \ftilde_n(\xb) \rvert > \epsilon)\\
	&= \PROB(\lvert \mutilde_n(A_n(\xb)) - \mu(A_n(\xb))  \rvert > \epsilon h_n^d)\\
	&= \PROB(\lvert G_{n,j}(1/2) - \EXP G_{n,j}(1/2) \rvert > \epsilon h_n^d (1 - 2H(-1/(2\sigma_W))))\\
	&= \PROB\left( \left\lvert \frac 1n \sum_{i=1}^n (\IND_{\{ W_{n,i,j \leq 1/2}\}} - \EXP [\IND_{\{ W_{n,i,j \leq 1/2}\}}]) \right\rvert > \epsilon h_n^d (1 - 2H(-1/(2\sigma_W)))\right) \\
	&\leq 2 \exp ( - \epsilon^2 n h_n^{2d} (1 - 2H(-1/(2\sigma_W)))^2/2).\end{align*}
The assumptions of the theorem now guarantee that
\begin{equation*}
	\sum_{n = 1}^\infty \PROB (\lvert \ftilde_n(\xb) - \EXP \ftilde_n(\xb) \rvert > \epsilon) < \infty,
\end{equation*}
and almost sure convergence follows using the Borel-Cantelli lemma.

We start the proof of $L^1$-consistency by recalling for $a \geq 0$ the identity $|a-b|=2(a-b \vee 0)_+ + 2b_- + b-a$ which is used to write
\begin{align*}
	&\int_{\R^d} |f(\xb) - \tilde f_n(\xb) | \dd \xb\\
	&\hspace{1em}=
	2\int_{\R^d} (f(\xb) - \tilde f_n(\xb) \vee 0)_+ \dd \xb
	+ 2 \int_{\R^d} (\tilde f_n(x))_- \dd \xb\\
	&\hspace{2em}+ \int_{\R^d} \tilde f_n(\xb)  \dd \xb -
	\int_{\R^d} f(\xb)  \dd \xb.
\end{align*}
The pointwise consistency together with Lebesgue's dominated convergence theorem implies that
\begin{align*}
	\int_{\R^d} (f(\xb) - \tilde f_n(\xb) \vee 0)_+ \dd \xb
	&\to 0 \quad \text{almost surely.}
\end{align*}
Thus, it remains to show that
\begin{align}\label{eq:proof:univ:cons:1}
	\int_{\R^d} \tilde f_n(x) \dd x
	&=
	\sum_{j=1}^{N_n} \tilde \mu_n(A_{h_n,j})
	\to 
	1 \quad \text{almost surely}
\end{align}
and
\begin{align}\label{eq:proof:univ:cons:2}
	\int_{\R^d} (\tilde f_n(x))_-  \dd x
	&\to 0 \quad \text{almost surely}.
\end{align}
We first show \eqref{eq:proof:univ:cons:1}.
By the triangle inequality, we have
\begin{align*}
	\left|\sum_{j=1}^{N_n} \tilde \mu_n(A_{h_n,j})-1\right|
	&\le 
	\left|\sum_{j=1}^{N_n}  \mu_n(A_{h_n,j})-1\right|
	+
	\left|\sum_{j=1}^{N_n} \tilde \mu_n(A_{h_n,j})
	-
	\sum_{j=1}^{N_n} \mu_n(A_{h_n,j})\right|.
\end{align*}
The first term on the right-hand side tends to $0$ almost surely (using the strong universal consistency of the standard histogram from Theorem~\ref{hcons}), and it remains to consider the second term.
Put
\begin{align*}
	\Delta_n &\defeq \left|\sum_{j=1}^{N_n} \tilde \mu_n(A_{h_n,j})-\sum_{j=1}^{N_n} \mu_n(A_{h_n,j})\right|.
\end{align*}
For any $\epsilon > 0$, consider the probability $\PROB(\Delta_n > \epsilon)$.
Conditioning on $\mathcal D_n$ yields
\begin{align*}
	\PROB(\Delta_n > \epsilon) = \EXP[\EXP[\IND_{\{ \Delta_n > \epsilon \}} | \mathcal D_n ]]. 
\end{align*}
Setting $C(\sigma_W) = (1 - 2H(-1/(2\sigma_W)))^{{-1}}$ implies
\begin{align*}
	\Delta_n &= C(\sigma_W) \left| \sum_{j=1}^{N_n} (G_{n,j}(1/2) - \EXP [G_{n,j}(1/2)|\mathcal D_n]) \right|\\
	&= C(\sigma_W)n^{-1} \left| \sum_{j=1}^{N_n} \sum_{i=1}^n (\IND_{\{ W_{n,i,j \leq 1/2}\}} - \EXP [\IND_{\{ W_{n,i,j \leq 1/2}\}}|\mathcal D_n]) \right|.
\end{align*}
Applying Hoeffding's inequality to the conditional distribution of the random variables $W_{n,i,j}$, $i=1,\ldots,n$, $j=1,\ldots,N_n$ given $\mathcal D_n$ (note that conditional on $\mathcal D_n$ all these random variables are independent) yields
\begin{align*}
	\EXP[\IND_{\{ \Delta_n > \epsilon \}} | \mathcal D_n ] &\leq 2 \exp \left( - \frac{2n\epsilon^2}{C(\sigma_W)^2 N_n} \right)
\end{align*}
Hence,
\begin{align*}
	\sum_{n=1}^\infty \PROB(\Delta_n > \epsilon) \lesssim \sum_{n=1}^\infty \exp \left( - \frac{2n\epsilon^2}{C(\sigma_W)^2 N_n} \right),
\end{align*}
and the sum on the right-hand side converges provided that $\frac{n}{N_n \log n} \to \infty$ which is equivalent to the condition $\frac{nh_n^d}{r_n^d \log n} \to \infty$ stated in the assumptions of the theorem (recall that $N_n \asymp r_n^d/h_n^d$).

In order to prove \eqref{eq:proof:univ:cons:2} it is sufficient to show $\sum_{n=1}^\infty \PROB(\lvert Z_n \rvert > \epsilon) < \infty$ for any $\epsilon > 0$ where $Z_n = \int_{\R^d} (\tilde f_n(x))_-  \dd x = \sum_{j=1}^{N_n} (\tilde \mu_n(A_{h_n,j}))_-$.
One has $\{ \lvert Z_n \rvert > \epsilon \} \subseteq \left\lbrace  \widetilde Z_n > \epsilon \right\rbrace$ where 
\begin{equation*}
	\widetilde Z_n \defeq C(\sigma_W) \sum_{j=1}^{N_n} \left[ H(1/(2\sigma_W)) - n^{-1} \sum_{i=1}^n \IND_{\{ \zeta_{i,j} \leq 1/(2\sigma_W)\} } \right]_-.
\end{equation*}
By a union bound argument and Hoeffding's inequality we obtain
\begin{align*}
	\sum_{n=1}^\infty& \PROB(\lvert Z_n \rvert > \epsilon) \leq \sum_{n=1}^\infty \PROB(\widetilde Z_n  > \epsilon)\\
	&\leq \sum_{n=1}^\infty N_n \PROB \left( H(1/(2\sigma_W)) - n^{-1} \sum_{i=1}^n \IND_{\{ \zeta_{i,1} \leq 1/(2\sigma_W)\} } < -\epsilon/(N_n C(\sigma_W)) \right)\\
	&\leq \sum_{n=1}^\infty N_n \exp \left( - \frac{2n\epsilon^2}{C^2(\sigma_W)N_n^2} \right),
\end{align*}
and the series converges under the additional assumption that $n/(N_n^2 \log n) \to \infty$ which is equivalent to $nh_n^{2d}/(r_n^{2d} \log n) \to \infty$.

\subsection{Proof of Theorem~\ref{Thm:prate}}

We begin with some preliminaries that will allow us to show that the upper bound obtained for the estimator $\hat f_n$ holds also for the estimator $\tilde f_n$.
For this purpose, note that
\begin{align*}
	\Var (\hat \mu_n(A_{h_n,j}))
	=\frac{\mu(A_{h_n,j})(1-\mu(A_{h_n,j}))+\sigma_W^2}{n}
\end{align*}
and
\begin{align*}
	\Var (\tilde \mu_n(A_{h_n,j}))
	&=\frac{\EXP[G_{n,j}(1/2)](1-\EXP[G_{n,j}(1/2)])}{n(1-2H(-1/(2\sigma_W)))^2}.
\end{align*}
Thus, the inequality
\begin{align}\label{eq:ineq:var}
	\Var (\tilde \mu_n(A_{h_n,j}))
	\le 
	\Var (\hat \mu_n(A_{h_n,j})) 
\end{align}
is equivalent to
\begin{align}
	\label{eg}
	\frac{\EXP[G_{n,j}(1/2)](1-\EXP[G_{n,j}(1/2)])}{(1-2H(-1/(2\sigma_W)))^2}
	&\le 
	\mu(A_{h_n,j})(1-\mu(A_{h_n,j}))+\sigma_W^2.
\end{align}
Note that
\begin{align*}
	\EXP[G_{n,j}(1/2)]
	&=
	\mu(A_{h_n,j})H(-1/(2\sigma_W))+(1-\mu(A_{h_n,j}))H(1/(2\sigma_W))\\
	&=
	\mu(A_{h_n,j})H(-1/(2\sigma_W))+(1-\mu(A_{h_n,j}))(1-H(-1/(2\sigma_W)))\\
	&=
	(2\mu(A_{h_n,j})-1)H(-1/(2\sigma_W))+1-\mu(A_{h_n,j})\\
	&=
	(1/2-\mu(A_{h_n,j}))(1-2H(-1/(2\sigma_W)))+1/2
\end{align*}
and
\begin{align*}
	H(-1/(2\sigma_W))
	&=
	\frac{1}{2}\exp(-1/(\sqrt{2}\sigma_W)).
\end{align*}
Therefore,
\begin{align*}
	&\frac{\EXP[G_{n,j}(1/2)](1-\EXP[G_{n,j}(1/2)])}{(1-2H(-1/(2\sigma_W)))^2}\\
	&=
	\frac{1/4-(1/2-\mu(A_{h_n,j}))^2(1-2H(-1/(2\sigma_W)))^2 }{(1-2H(-1/(2\sigma_W)))^2}\\
	&=
	\frac{1/4}{(1-2H(-1/(2\sigma_W)))^2}-(1/2-\mu(A_{h_n,j}))^2,
\end{align*}
and so \eqref{eg} holds if and only if
\begin{align*}
	\frac{1/4}{(1-2H(-1/(2\sigma_W)))^2}
	&\le 
	(1/2-\mu(A_{h_n,j}))^2+\mu(A_{h_n,j})(1-\mu(A_{h_n,j}))+\sigma_W^2,
\end{align*}
or equivalently,
\begin{align}\label{eq:inequality:sigma_W}
	\frac 14\left(\frac{1}{\left(1-\exp\left(-1/(\sqrt{2}\sigma_W)\right)\right)^2}-1\right)
	&\le 
	\sigma_W^2.
\end{align}
This last inequality (which does not involve $\mu(A_{h_n,j})$ anymore) holds for any choice of $\sigma_W$ which can be seen as follows: the ratio of the left-hand side and the right-hand side of \eqref{eq:inequality:sigma_W} converges to $0$ for $\sigma_W \to 0$, to $1/2$ for $\sigma_W \to \infty$, and has a unique maximum on $(0,\infty)$ at $\sigma_W \approx 0.614$ with value $0.755 \leq 1$.
Hence, we have proven \eqref{eq:ineq:var} which together with the bound of the bias term in Theorem~\ref{hrate} implies
\begin{align}
	&\EXP \left[  \int_{\R^d} |f(\xb ) - \ftilde_n(\xb) | \dd \xb\right] \notag\\ 
	&\leq 
	\EXP\left[  \int_{\R^d} |f(\xb ) - \EXP\{\ftilde_n(\xb)\} | \dd \xb\right]  
	+
	\EXP\left[  \int_{\R^d} |\EXP\{\ftilde_n(\xb)\} - \ftilde_n(\xb) | \dd \xb\right] \notag\\
	&\leq 
	\EXP\left[  \int_{\R^d} |f(\xb ) - \EXP\{f_n(\xb)\} | \dd \xb\right]  
	+\int_{\R^d \setminus \bigcup_{j=1}^{N_n} A_{h_n,j}} f(\xb) \dd \xb\notag\\
	&\qquad+
	\sum_{j=1}^{N_n} \sqrt{\Var (\tilde \mu_n(A_{h_n,j}))}\notag\\
	&\leq 
	L\sqrt{d}h_n 
	+\int_{\R^d \setminus \bigcup_{j=1}^{N_n} A_{h_n,j}} f(\xb) \dd \xb
	+
	\sum_{j=1}^{N_n} \sqrt{\Var (\hat \mu_n(A_{h_n,j}))}\notag\\
	&\leq 
	L\sqrt{d}h_n 
	+\int_{\R^d \setminus \bigcup_{j=1}^{N_n} A_{h_n,j}} f(\xb) \dd \xb
	+
	\sum_{j=1}^{N_n} \frac{\sqrt{\mu(A_{h_n,j})+\sigma_W^2}}{\sqrt{n}}.\notag \end{align}
Note that the assumptions of the theorem guarantee the existence of a constant $C_L$ (depending on the Lipschitz constant $L$) such that $\lVert f \rVert_\infty \leq C_L$ holds.
In combination with $\sqrt{a+b} \leq \sqrt a + \sqrt b$, this implies
\begin{align*}
	&\EXP\left[  \int_{\R^d} |f(\xb ) - \ftilde_n(\xb) | \dd \xb\right] \\
	&\leq L\sqrt{d}h_n  +\int_{\R^d \setminus \bigcup_{j=1}^{N_n} A_{h_n,j}} f(\xb) \dd \xb + \frac{N_n\sqrt{C_L h_n^d}}{\sqrt{n}} + \frac{N_n \sigma_W}{\sqrt{n}}\\
	&\leq L\sqrt{d}h_n + \int_{\R^d \setminus \bigcup_{j=1}^{N_n} A_{h_n,j}} f(\xb) \dd \xb + O \left( \frac{r_n^d}{\sqrt{nh_n^d}} \right) + O \left( \frac{\sigma_Wr_n^d}{\sqrt{nh_n^{2d}}} \right).
\end{align*}
Exploiting that $\EXP[\ftilde_n(\xb)] = \EXP[\hat f_n(\xb)]$, one can use the same argument in order to show that the same upper bound holds for the estimator $\hat f_n$.
Now one can see directly that we can attain the classical rate if $\sigma_W$ is sufficiently small.

\subsection{Proof of Theorem~\ref{THM:LOWER:LIPSCHITZ}}\label{proof:lower}
As usual in proving minimax lower bounds, we bound the complexity of the whole estimation problem by restricting the class of potential density functions to a finite subset which is parametrised by a vector $\theta \in \{ \pm 1\}^{k^d}$ for some $k$ that will be specified later. It is well-known in nonparametric statistics that for proving global lower bounds the reduction to two hypothesis (which works for deriving sharp lower bounds for the estimation of linear functionals of the density) is not sufficient. The main additional ingredient in contrast to classical lower bounds is the data processing inequality \eqref{eq:data:processing:proof} that explains the deterioration of the rate of convergence.

In the sequel, whenever $\mu_\theta$ denotes the distribution of the original i.i.d.\,data $X_1,\ldots,X_n$ with density $f_\theta$, we denote for a privacy mechanism $Q$ with $\mu_\theta^{\Zb}$ the resulting distribution of the privatised data $\Zb = (Z_1,\ldots,Z_n)$ (the dependence of this distribution on the mechanism $Q$ is suppressed for the sake of convenience).

Let $\PROB$ and $\mathbb Q$ be two probability measures on the same measurable space $(\mathcal X, \mathscr X)$, and denote by $p$ and $q$ their densities with respect to some dominating measure $\nu$.
Recall that the total variation distance, the Hellinger distance, and the Kullback-Leibler distance are defined as
\begin{align*}
	&\TV(\PROB, \mathbb Q) = \frac{1}{2} \int_{\mathcal X} \lvert p(
	\xb) - q(\xb) \rvert \nu(\dd \xb),\\
	&H(\PROB, \mathbb Q) = \left( \int_{\mathcal X} ( \sqrt{p(\xb)} - \sqrt{q(\xb)} )^2 \nu(\dd \xb) \right)^{1/2},\\
	&\KL(\PROB, \mathbb Q) = \int_{\mathcal X} \log \left( \frac{p(\xb)}{q(\xb)} \right) p(\xb) \nu(\dd \xb),
\end{align*}
respectively (the last definition being valid for $\PROB \ll \mathbb Q$; otherwise one defines $\KL(\PROB, \mathbb Q)= + \infty$).
For the proof we need the following result which links the Kullback-Leibler distance of the privatised observations $Z_1,\ldots,Z_n$ (generated by a potentially interactive privacy mechanism) to the total-variation distance of the original random variables.

\begin{lemma}[Consequence of \citep{DuJoWa18}, Theorem~1]\label{LEM:LOWER:EX:DUCHI} Let $\alpha \geq 0$.
	For any $\alpha$-locally differentially private mechanism $Q \in \Qc_\alpha$
	\begin{equation*}
		\KL(\mu_\theta^{\Zb},\mu_{\theta^\prime}^{\Zb}) \leq 4n(e^\alpha -1)^2  \TV^2(\mu_\theta, \mu_{\theta^\prime}).
	\end{equation*}
\end{lemma}

For the proof of Theorem~\ref{THM:LOWER:LIPSCHITZ}, we need to introduce some notation first.
Define the function $g_0 \colon [0,1]^d \to \R$ by
\begin{equation*}
	g_0(\xb) = L \min_{i=1,\ldots,d} \min \{ x_i, 1-x_i \}.
\end{equation*}
Then, denoting with $h$ the $1$-Lipschitz function $h\colon [0,1] \to \R$ defined by $h(x) = x \IND_{ \{ x \in [0,1/2] \}} + (1-x)\IND_{\{x \in (1/2,1]\} }$, we have
\begin{align*}
	\lvert g_0(\xb) - g_0(\xb\p)\rvert &= L \,\lvert \min_{i=1,\ldots,d} \min \{ x_i, 1-x_i \} - \min_{i=1,\ldots,d} \min \{ x_i\p, 1-x_i\p \} \rvert\\
	&\leq L \max_{i=1,\ldots,d} \, \lvert \min \{ x_i, 1-x_i \} - \min \{ x_i\p, 1-x_i\p \} \rvert \\
	&= L \max_{i=1,\ldots,d}\, \lvert h(x_i) - h(x_i\p) \rvert\\
	&\leq L \max_{i=1,\ldots,d} \lvert x_i - x_i\p \rvert\\
	&\leq L \lVert \xb - \xb\p \rVert,
\end{align*}
showing that $g_0$ belongs to $\Flipd(L)$.

For a positive integer $k$ (that will be specified below) set $A=[0,1/(2k))^d$. Define $\yb_\jb = (y_{j_1},\ldots,y_{j_d}) \in [0,1]^d$ for $\jb \in \{ 0,\ldots,k-1 \}^d$ by $y_{j_i} = 1/4 + j_i/(2k)$, and further put $A_\jb = \yb_\jb + A$ for the same values of $\jb$.
Now, consider the function $g\colon A \to \R$ defined as
\begin{equation*}
	g(\xb) =  \frac{(-1)^{\# \{ i \in \{ 1,\ldots,d \} \, : \, 1/(4k) \leq x_i < 1/(2k) \} }}{4k} \cdot g_0(\tau(x_1),\ldots,\tau(x_d))
\end{equation*}
where
\begin{align*}
	\tau(x) = \begin{cases}
		4kx,& \text{ if } 0 \leq x < \frac{1}{4k},\\
		4k(x - 1/(4k)),& \text{ if } \frac{1}{4k} \leq x < \frac{1}{2k}.
	\end{cases}
\end{align*}
\begin{figure}[t]
	\centering
	\includegraphics[width=9cm]{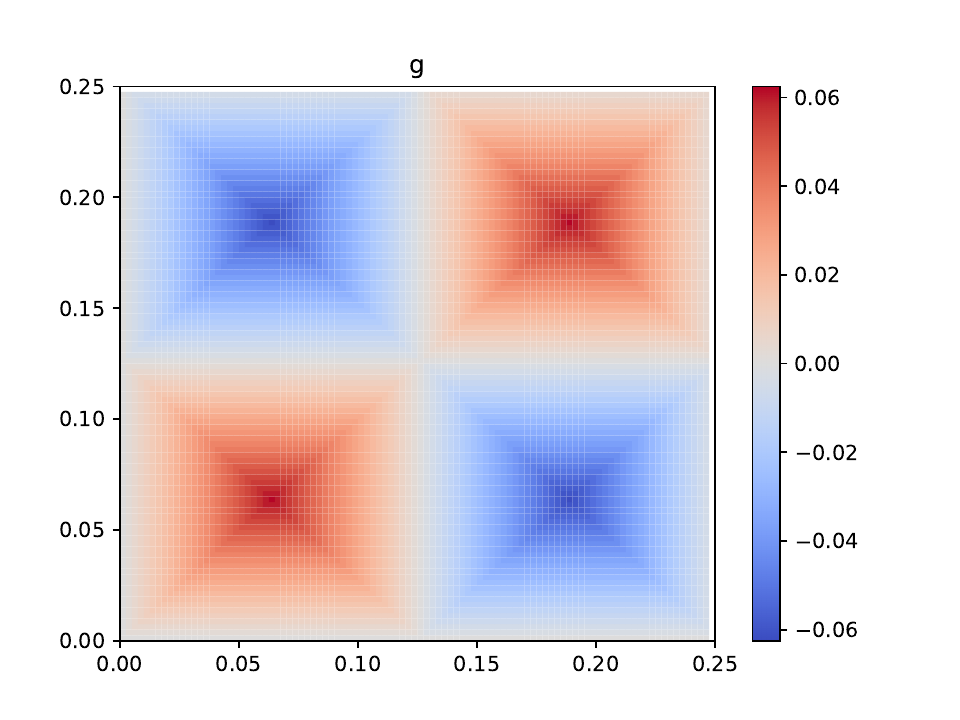}
	\caption{Heatmap of the function $g\colon [0,0.25]^2 \to \R$ for the case $d=k=2$ with $L=1$. This function is used to construct the hypotheses $f_\theta$ in the proof of the lower bound. Note that positive and negative parts of the functions cancel each other when integrating over the whole area $[0,0.25)^2$.}\label{fig:g}
\end{figure}
Take $\theta = (\theta_{\jb}) \in \{ \pm 1 \}^{k^d}$.
With such a $\theta$ we associate the function defined by
\begin{equation*}
	f_\theta (\xb) = \begin{cases}
		f_0(\xb), & \text{ if }\xb \notin \bigcup_{\jb} A_{\jb},\\
		f_0(\xb) + \theta_{\jb} g(\xb - \yb_\jb), &\text{ if } \xb \in A_{\jb}.
	\end{cases}
\end{equation*}
Here, the function $f_0$ is chosen such that it is constant on the block $\bigcup_{\jb} A_\jb = [1/4,3/4)^d$
, and outside this block it is defined in such a way that $f_0$ is positive, integrates to $1$ and satisfies the $\Flipd(L)$ condition.
The existence of such a function $f_0$ is guaranteed whenever $L$ is sufficiently large.
The definition of $g$ (see Figure~\ref{fig:g} for an illustration in the case $d=k=2$) guarantees that the hypotheses $f_\theta$ are density functions and belong to $\Flipd$.

After these preparations, we now consider privatised data $\Zb = (Z_1,\ldots,Z_n)$ generated by an arbitrary $\alpha$-LDP privacy mechanism $Q \in \Qc_\alpha$ (generating observations in an arbitrary measurable space $(\mathcal Z, \mathscr Z)$) and let us denote with $\ftilde$ an arbitrary estimator based on these data.
Then, denoting with $\EXP_f$ and $\EXP_\theta$ the expectation operator when the true density is $f$ and $f_\theta$, respectively, we have
\begin{align}
	\sup_{f \in \Flipd} \EXP_f \left[ \int_{[0,1]^d}
	\lvert \ftilde(\xb) - f(\xb) \rvert \dd\xb \right] &\geq \sup_{\theta \in \{ \pm 1 \}^{k^d}} \EXP_\theta \left[ \int_{[0,1]^d} \lvert \ftilde(\xb) - f_\theta(\xb) \rvert \dd\xb \right]\notag\\
	&\hspace{-4em}\geq \frac{1}{2^{k^d}} \sum_{\theta \in \{ \pm 1 \}^{k^d}} \EXP_\theta \left[ \int_{[0,1]^d} \lvert \ftilde(\xb) - f_\theta(\xb) \rvert \dd\xb \right]\notag\\
	&\hspace{-4em}\geq \frac{1}{2^{k^d}} \sum_{\theta \in \{ \pm 1 \}^{k^d}} \sum_{\jb} \EXP_\theta \left[ \int_{A_\jb} \lvert \ftilde(\xb) - f_\theta(\xb) \rvert \dd\xb \right]\notag\\
	&\hspace{-4em}= \frac{1}{2^{k^d}} \sum_{\jb} \sum_{\theta \in \{ \pm 1 \}^{k^d}} \EXP_\theta \left[ \int_{A_\jb} \lvert \ftilde(\xb) - f_\theta(\xb) \rvert \dd\xb \right]\notag\\
	&\hspace{-4em}= \frac{1}{2^{k^d+1}} \sum_{\jb} \sum_{\theta \in \{ \pm 1 \}^{k^d}} \left\lbrace \EXP_\theta \left[ \int_{A_\jb} \lvert \ftilde(\xb) - f_\theta(x) \rvert \dd\xb \right] \right.\notag\\
	&\hspace{-4em}\left. + \, \EXP_{\theta^{(\jb)}} \left[ \int_{A_\jb} \lvert \ftilde(\xb) - f_{\theta^{(\jb)}}(\xb) \rvert \dd\xb \right] \right\rbrace, \label{eq:sup:red}
\end{align}
where we denote for some given $\theta \in \{\pm 1 \}^{k^d}$ with $\theta^{(\jb)}$ the element of $\{\pm 1 \}^{k^d}$ satisfying $\theta^{(\jb)}_\ib = \theta_\ib$ for $\ib \neq \jb$ and $\theta^{(\jb)}_\jb = -\theta_\jb$ (sign reversal at the coordinate with index $\jb$).
Consider the Hellinger affinity 
\[
\rho(\mu_{\theta}^{\Zb}, \mu_{\theta^{(\jb)}}^{\Zb}) \defeq \int_{\mathcal Z^n} \sqrt{ \psi_\theta^{\Zb}(\zb) \psi_{\theta^{(\jb)}}^{\Zb} (\zb)} \nu(\dd \zb)
\]
where $\psi_\theta^{\Zb}$ and $\psi_{\theta^{(\jb)}}^{\Zb}$ denote densities of the measure $\mu_{\theta}^{\Zb}$ and $\mu_{\theta^{(\jb)}}^{\Zb}$, respectively, with respect to some dominating measure, say $\nu$.
Using the elementary inequality $\sqrt{a+b} \leq \sqrt a + \sqrt b$ for $a,b \geq 0$ and the Cauchy-Schwarz inequality, we obtain 
\begin{align*}
	\rho(\mu_{\theta}^{\Zb}, \mu_{\theta^{(\jb)}}^{\Zb}) &= \int_{\mathcal Z^n} \frac{\left( \int_{A_\jb} \lvert f_\theta(\xb) - f_{\theta^{(\jb)}}(\xb) \rvert \dd \xb \right)^{1/2}}{\left( \int_{A_\jb} \lvert f_\theta(\xb) - f_{\theta^{(\jb)}}(\xb) \rvert \dd \xb \right)^{1/2}} \sqrt{\psi_\theta^{\Zb}(\zb) \psi_{\theta^{(\jb)}}^{\Zb}(\zb)} \nu(\dd \zb)\\
	&\leq \int_{\mathcal Z^n} \frac{\left( \int_{A_\jb} \lvert f_\theta(\xb) - \ftilde(\xb) \rvert \dd \xb \right)^{1/2}}{\left( \int_{A_\jb} \lvert f_\theta(\xb) - f_{\theta^{(\jb)}}(\xb) \rvert \dd\xb \right)^{1/2}} \sqrt{\psi_\theta^{\Zb}(\zb)}\sqrt{\psi_{\theta^{(\jb)}}^{\Zb}(\zb)} \nu(\dd \zb)\\
	&\hspace{1em}+ \int_{\mathcal Z^n} \frac{\left( \int_{A_\jb} \lvert f_{\theta^{(\jb)}}(\xb) - \ftilde(\xb) \rvert \dd \xb \right)^{1/2}}{\left( \int_{A_\jb} \lvert f_\theta(\xb) - f_{\theta^{(\jb)}}(\xb) \rvert \dd \xb \right)^{1/2}}  \sqrt{\psi_{\theta^{(\jb)}}^{\Zb}(\zb)} \sqrt{\psi_\theta^{\Zb}(\zb)}\nu(\dd \zb)\\
	&\leq \left( \int_{\mathcal Z^n} \frac{\int_{A_\jb} \lvert f_\theta(\xb) - \ftilde(\xb) \rvert \dd \xb}{ \int_{A_\jb} \lvert f_\theta(\xb) - f_{\theta^{(\jb)}}(\xb) \rvert \dd \xb} \psi_\theta^{\Zb}(\zb)\nu(\dd \zb) \right)^{1/2}\\
	&\hspace{1em}+ \left( \int_{\mathcal Z^n} \frac{\int_{A_\jb} \lvert \ftilde(\xb) - f_{\theta^{(\jb)}}(\xb) \rvert \dd \xb }{ \int_{A_\jb} \lvert f_\theta(\xb) - f_{\theta^{(\jb)}}(\xb) \rvert \dd \xb } \psi_{\theta^{(\jb)}}^{\Zb}(\zb)\nu(\dd \zb)\right)^{1/2}.
\end{align*}
Consequently, using $(a+b)^2 \leq 2a^2 +2b^2$,
\begin{align}
	\frac 1 2\rho^2(\mu_{\theta}^{\Zb}, \mu_{\theta^{(j)}}^{\Zb}) \int_{A_\jb} \lvert f_\theta(\xb) - f_{\theta^{(\jb)}}(\xb) \rvert \dd\xb\notag\\
	&\hspace{-15em}\leq \EXP_\theta \left[ \int_{A_\jb} \lvert f_\theta(\xb) - \ftilde(\xb) \rvert \dd\xb \right] + \EXP_{\theta^{(\jb)}} \left[ \int_{A_\jb} \lvert f_{\theta^{(\jb)}}(\xb) - \ftilde(\xb) \rvert \dd\xb \right].\label{eq:affinity:EXP}
\end{align}
Let us now determine a lower bound for the quantity $\rho(\mu_{\theta}^{\Zb}, \mu_{\theta^{(\jb)}}^{\Zb})$.
For this, we first derive an upper bound for the Kullback-Leibler divergence $\KL (\mu_{\theta}^{\Zb},\mu_{\theta^{(\jb)}}^{\Zb})$.
Lemma~\ref{LEM:LOWER:EX:DUCHI} directly provides the estimate
\begin{equation}\label{eq:data:processing:proof}
	\KL(\mu_{\theta}^{\Zb}, \mu_{\theta^{(j)}}^{\Zb}) \leq 4n(e^\alpha - 1)^2 \TV^2(\mu_{\theta},\mu_{\theta^{(j)}}).
\end{equation}
Further, by definition of $f_\theta$ and $f_{\theta^{(\jb)}}$,
\begin{align*}
	\TV(\mu_{\theta},\mu_{\theta^{(\jb)}}) &= \frac{1}{2} \int_{A_\jb} \lvert f_\theta(\xb) - f_{\theta^{(\jb)}}(\xb) \rvert \dd \xb\\
	&= \int_A \lvert g(\xb) \rvert \dd \xb\\
	&= \frac{2^d}{4k} \int_{[0,1/(4k))^d} g_0(4kx_1,\ldots,4kx_d) \dd x_1 \ldots \dd x_d\\
	&= \frac{2^d}{(4k)^{d+1}} \int_{[0,1)^d} g_0(\xb) \dd \xb.
\end{align*}
Now choose $k$ as the smallest integer $\geq 1$ such that
\begin{equation*}
	(n(e^\alpha -1)^2 \vee 1)\frac{1}{(2k)^{2d+2}} \left( \int_{[0,1)^d} g_0(\xb) \dd \xb  \right)^2 \leq 1,
\end{equation*}
which yields $\KL(\mu_{\theta}^{\Zb}, \mu_{\theta^{(j)}}^{\Zb}) \leq 1$.
Using the relation $H^2 \leq \KL$ between squared Hellinger and Kullback-Leibler distance (cf.\,\cite{tsybakov2009introduction}, Equation~(2.19)) together with the identity $\rho = 1- \frac{1}{2}H^2$ for the Hellinger affinity (cf.\,\cite{tsybakov2009introduction}, p.~83, property~(iii) of the Hellinger distance), we obtain the estimate
\begin{equation*}
	\rho(\mu_{\theta}^{\Zb}, \mu_{\theta^{(\jb)}}^{\Zb}) \geq \frac{1}{2},
\end{equation*}
which is independent of the index $\jb$.
Combining this last estimate with \eqref{eq:affinity:EXP} and putting the result into \eqref{eq:sup:red} yields
\begin{align*}
	\sup_{f \in \Flipd(L)} \EXP_f \left[ \int_{[0,1]^d} \lvert \ftilde(\xb) - f(\xb) \rvert \dd \xb \right] &\geq \frac{1}{8} \sum_{\jb} \int_{A_\jb} \lvert f_\theta(\xb) - f_{\theta^{(\jb)}}(\xb) \rvert \dd \xb\\
	&= \frac{k^d}{4} \int_{A} \lvert g(\xb) \rvert \dd \xb\\
	&\asymp k^{-1}\\
	&\asymp (n(e^\alpha - 1)^2)^{-\frac{1}{2d+2}} \wedge 1,
\end{align*}
where we used that $k \asymp (n(e^\alpha - 1)^2 \vee 1)^{\frac{1}{2d+2}}$.
This proves the desired lower bound.

\end{document}